\documentclass[12pt,a4paper]{article}
\usepackage{amsmath}
\usepackage{amsthm}
\usepackage{amsfonts}
\usepackage{amssymb}
\usepackage{stmaryrd}
\usepackage{latexsym}
\usepackage[usenames]{color}
\usepackage[russian, english]{babel}

\theoremstyle{definition}

\newtheorem{example}{Example}[section] 

\begin{document}

\begin{center}
\Large{
An analogue of the ElGamal scheme based on the Markovski algorithm}

\large{
Nadeghda Malyutina, Victor Shcherbacov}

\end{center}
\begin{abstract}
We give an analogue of the ElGamal encryption system based on the Markovski algorithm [4; 5].

AMS Classification 20N05, 05B15, 94A60


\end{abstract}


\hfill{In memory of academician M. M. Choban}

\section {Introduction}

Usually the classical Taher ElGamal encryption system is formulated in the language of number theory using multiplication modulo a prime [1].

ElGamal’s scheme is a public key cryptosystem based on the difficulty of computing discrete logarithms in a finite field. The cryptosystem includes an encryption algorithm and a digital signature algorithm. ElGamal Scheme Underlies US Former Electronic Digital Signature Standards (DSA) and Russia (GOST R 34.10-94).

The scheme was proposed by Taher ElGamal in 1985. ElGamal developed one of the variants of the Diffie-Hellman algorithm. He improved the Diffie-Hellman system and obtained two algorithms that were used for encryption and for authentication. Unlike RSA, the ElGamal algorithm was not patented and, therefore, became a cheaper alternative, since it did not require payment of license fees [2].

The sender of messages and their recipient can be individuals, organizations, or technical systems. These may be subscribers of a network, users of a computer system, or abstract “parties” involved in information interaction. But more often participants are identified with people and replaced with the formal designations A and B by Alice and Bob.
It is assumed that messages are transmitted through the so-called “open” communication channel, available for listening to some other persons.

In cryptography, it is usually assumed that a person sending messages or receiving them has some opponent E and this opponent can intercept messages transmitted over an open channel. The enemy is considered as a certain person named Eve, who has at her disposal powerful computing equipment and owns cryptanalysis methods. Naturally, Alice and Bob want their messages to be incomprehensible to Eve, and use special ciphers for this.

Before sending a message over an open communication channel from A to B, A encrypts the message, and B, having received the encrypted message, decrypts it, restoring the original text. The important thing is that Alice and Bob can agree on the cipher they use not on an open channel, but on a special "closed" channel, inaccessible for listening to the enemy. It should be borne in mind that usually the organization of such a closed channel and the transmission of messages through it is too expensive compared to an open channel or a closed channel cannot be used at any time.
Each attempt to break the cipher is called an attack on the cipher. In cryptography, it is generally accepted that the adversary can know the encryption algorithm used, the nature of the transmitted messages and the intercepted ciphertext, but does not know the secret key.

Developers of modern cryptosystems strive to make attacks on known and selected text invulnerable. Significant successes have been achieved along this path.

\section {ElGamal's scheme}

Suppose there are subscribers $A, B, C, \dots$  who want to transmit encrypted messages to each other without having any secure communication channels. We will consider the code proposed by ElGamal, which solves this problem, using, in contrast to the Shamir code, only one message forwarding. In fact, the Diffie -Hellman scheme is used here to form a common secret key for two subscribers transmitting a message to each other, and then the message is encrypted by multiplying it by this key. For each subsequent message, the secret key is recalculated.
A large prime number is selected $p$ and number $g$, such that different degrees of $g$ are different numbers modulo p. The numbers p and g are transmitted to subscribers in the clear.

Then each subscriber of the group selects his secret number $c_i$,    $1 < c_i< p-1$, and calculates the corresponding open number $d_i$,

\begin{equation}
d_i=g^{c_i} \pmod p \label{1_{Ro}}
\end{equation}

\begin{center}
Table 1.  (User keys in the ElGamal system)
$$
\begin{array}{|l|l|l|}
\hline
$\text Subscriber $	&  $\text Private \, key $	& $\text Public \: key $ \\
\hline
A	& c_A	& d_A \\
\hline
B	& c_B &	d_B \\
\hline
C &	c_C &	d_C \\
\hline
\end{array}
$$
\end{center}

We show now how A sends message m to the subscriber B.
We will assume, that the message is presented as a number $m < p$.

Step 1.  A forms a random number $k, 1 \geq k \geq p-2$, calculates numbers:
\begin{equation}
r=g^k  \mod p \label{2_{Ro}}
\end{equation}

\begin{equation}
e=m \cdot {d_{B}}^k  \mod p	\label{3_{Ro}}		
\end{equation}

and passes a couple of numbers $(r,e)$ to the subscriber $B$.

Step 2  B, getting  $(r,e)$, calculates
\begin{equation}
m^{\prime}=e\cdot r^{p-1-c_B}  \mod p  	
\end{equation}

Statement 1 (properties of the ElGamal cipher).

(1)	The subscriber B received a message, i.e. $m^{\prime} = m$;

(2)	the adversary, knowing $p,g,d_B,r$ and $e$, cannot calculate $m$.

\begin{example}  Consider the transmission of message m = 15 from A to B.

We take  p = 23,  g = 5. Let subscriber B choose for himself a secret number $c_B=13$ and calculate (\ref{1_{Ro}}):
$d_B=5^{13} \pmod {23} =21$.

Subscriber A randomly selects the number $k$, for example, $k = 7$, and calculates from (\ref{2_{Ro}}), (\ref{3_{Ro}}): $r = 5^7 \pmod {23} = 17$, $e = 15 \cdot  21^7  \pmod {23}= 15 \cdot 10 \pmod {23} = 12.$

Now $A$  sends to $B$  an encrypted message in the form of a pair of numbers $(17,12)$ and $B$ calculates:
$m^{\prime}=12 \cdot 17^{\,23-1-13}  \pmod {23} =    12\cdot 17^{9} \pmod{ 23}=12\cdot 7 \pmod {23}=15.$
So $B$ was able to decrypt the transmitted message.
\end{example}

By a similar scheme, all subscribers in the network can send messages. Moreover, any subscriber who knows the public key of subscriber B can send him messages encrypted using the public key $d_B$.
But only subscriber B, and no one else, can decrypt these messages using the secret key $c_B$ known only to him.

The Shamir cipher completely solves the problem of exchanging messages that are closed for reading, in the case when subscribers can use only open communication lines.

However, this message is sent three times from one subscriber to another, which is a drawback. The ElGamal cipher allows you to solve the same problem in one data transfer, but the amount of transmitted ciphertext is twice the size of the message.

It is easy to see that this system can also be formulated in terms of a residue ring modulo $p$ or, equivalently, using the language of the Galois field  $GF(p)$.

In addition, we can use the concept of the action of a group of automorphisms of a cyclic group $(Z_p,+)$ on this group.
Let $(Z_p,+)$ be a cyclic group of residues of large  simple order with respect to addition of residues and element $a$ be the generator of the group $(Z_{(p-1)},\cdot)\cong Aut(Z_p,+)(gcd(a,p-1)=1)$.

Alice's keys are the following: Public key $p$, $a$ and $a^m$, $m\in \mathbb{N}$.
Private key $m$.

Encryption. To send a message $b \in  (Z_{(p-1)},\cdot )$, Bob is calculating $a^r$ and  $a^{mr}$ for random $r\in \mathbb N$ (sometimes the number $r$ is called an ephemeral key [3]).

Ciphertext: $(a^r; a^{mr} \cdot b)$.

Decryption. Alice knows $m$, so if she gets the ciphertext $(a^r; a^{mr} \cdot  b)$, she will calculate $a^{mr}$ from $a^r$ and then $a^{(-mr)}$ and then from $a^{mr} \cdot b$ calculate {b}.

\begin{example}
Example  Alice picks $p=107, a= 2, m=67,$ and calculates
$a^m= 2^{67} \equiv  94 \pmod {107}$.

Her public key $(p,a^m)=(107,94)$, and her private key is $m=67$.
Bob wants to send a message \lq\lq B \rq\rq  to Alice. He selects a random integer $r=45$ and encrypts $B= 66$ like  $(a^m)^r  \cdot B$.

Bob gets:
$(2^{45},94^{45} \cdot 66) \equiv (28, 5 \cdot 66) \equiv (28, 9) \pmod {107}$.

He sends an encrypted message $(28,9)$ to Alice. Alice receives this message and using her private key $m = 67$ she decrypts as follows:
$28^{(-67)}\cdot  9 = 28^{(106-67)} \cdot 9 \equiv 28^{39}  \cdot 9 \equiv 43\cdot 9 \equiv 66 \pmod {107}.$
\end{example}

The complexity of this system is based on the complexity of the discrete logarithm problem. ElGamal's encryption system is not secure according to the selected attack ciphertext [3]. ElGamal cryptosystems are usually used in a hybrid cryptosystem, i.e. the message itself is encrypted using a symmetric cryptosystem and ElGamal also uses a symmetric cryptosystem to encrypt the key.

\section{An analogue of the ElGamal scheme based on the Markovski algorithm}

We give an analogue of the ElGamal encryption system based on the Markovski algorithm [4; 5].

Let $(Q,f)$ be a binary quasigroup and $T = (\alpha, \beta, \gamma)$ its isotopy.

Alice's keys are as follows:
The public key is $(Q,f)$, $T$, $T^{(m,n,k)} = (\alpha^m,\beta^n,\gamma^k)$ , $m, n, k \in \mathbb N$, and the Markovski algorithm.

Private key $m, n, k$.

Encryption.
To send a message $b \in (Q,f)$, Bob  calculated  $T^{(r,s,t) }$, $T^{(mr,ns,kt)}$ for random $r,s,t \in \mathbb N$ and $(T^{(mr,ns,kt)} (Q,f))$.

The ciphertext is  $(T^{(r,s,t)} ,T^{(mr,ns,kt)} (Q,f)b).$

To obtain $(T^{(mr,ns,kt)}  (Q,f)b)$, Bob uses the Markovski algorithm which is known to Alice.

Decryption
Alice knows $m,n,k,$  so if she gets the ciphertext $(T^{(r,s,t) },(T^{(mr,ns,kt)} (Q, f))b)$, she will calculate  $(T^{(mr,ns,kt)} (Q,f))^{(-1)}$  using  $T^{(r,s,t)}$  and finally she will calculate $b$.

\begin{example}  Let $(Q, f)$ be a binary quasigroup defined by the following Cayley table:

\smallskip

\centerline{Table 2.}
$$
\begin{array}{|l|l|l|l|l|l|l|l|}
\hline
f &	0	&1	&2	&3	&4	&5	&6 \\
\hline
0	&5	&2	&6	&4	&0	&3	&1 \\
\hline
1	&1	&6	&5	&3	&4	&2	&0 \\
\hline
2	&0	&5	&4	&6	&3	&1	&2 \\
\hline
3	&4	&1	&3	&0	&2	&6	&5 \\
\hline
4	&2	&4	&0	&1	&6	&5	&3 \\
\hline
5	&6	&3	&1	&2	&5	&0	&4 \\
\hline
6	&3	&0	&2	&5	&1	&4	&6\\
\hline
\end{array}
$$
and $T=(\alpha, \beta, \gamma)$ its isotopy, where:
$\alpha =(2 3 4)  (0 5 1 6)$  corresponds to a permutation of rows of a quasigroup table $Q$;
$\beta=(0 3 2 1)  (5 6)$  corresponds to a permutation of the columns of a quasigroup table $Q$ obtained after application $\alpha$;
$\gamma=(1 2 3 6 0 5 4)$   substitution applied to the table obtained after application  $\beta$.

And for $\gamma$ we have the inverse $\gamma^{(-1)}$  the following kind:
$\gamma^{(-1)} = (1 4 5 0 6 3 2) $.

Cayley tables of these permutations are of the form:

\centerline{Table 3.}
$$
\begin{array}{|l|l|l|l|l|l|l|l|}
\hline
\alpha 	&0	&1	&2	&3	&4	&5	&6\\
\hline
0	&6	&3	&1	&2	&5	&0	&4 \\
1	&3	&0	&2	&5	&1	&4	&6 \\
2	&4	&1	&3	&0	&2	&6	&5 \\
3	&2	&4	&0	&1	&6	&5	&3 \\
4	&0	&5	&4	&6	&3	&1	&2 \\
5	&1	&6	&5	&3	&4	&2	&0 \\
6	&5	&2	&6	&4	&0	&3	&1 \\
\hline
\end{array}
$$

\centerline{Table 4.}
$$
\begin{array}{|l|l|l|l|l|l|l|l|}
\hline
\beta	&0	&1	&2	&3	&4	&5	&6 \\
\hline
0	&2	&6	&3	&1	&5&  4&	0  \\
1	&5	&3	&0	&2	&1&	6&	4  \\
2	&0	&4	&1	&3	&2&	5&	6  \\
3	&1	&2	&4	&0	&6&	3&	5  \\
4	&6	&0	&5	&4	&3&	2&	1  \\
5	&3	&1	&6	&5	&4&	0&	2  \\
6	&4	&5	&2	&6	&0&	1&	3  \\
\hline
\end{array}
$$

\centerline{Table 5.}
$$
\begin{array}{|l|l|l|l|l|l|l|l|}
\hline
\gamma^{-1} &0	&1	&2	&3	&4	&5	&6 \\
\hline
0&	1&	3&	2&	4&	0&	5&	6    \\
1&	0&	2&	6&	1&	4&	3&	5    \\
2&	6&	5&	4&	2&	1&	0&	3    \\
3&	4&	1&	5&	6&	3&	2&	0    \\
4&	3&	6&	0&	5&	2&	1&	4    \\
5&	2&	4&	3&	0&	5&	6&	1    \\
6&	5&	0&	1&	3&	6&	4&	2    \\
\hline
\end{array}
$$

Then Alice’s keys are as follows:
The private key:  $m=3, n=6, k=5$.
The public key is $(Q,f),T,T^{(3,6,5)}=(\alpha^3,\beta^6,\gamma^5)$ and the Markovski algorithm, where:
$\alpha^3=(0 6 1 5);   \beta^6=(0 2)  (1 3);   \gamma^5=(0 3 1 5 6 2 4),
 \gamma^{-5}=(0 4 2 6 5 1 3)$.

 As a result, we get the following Cayley tables:
  \newpage
 \centerline{Table 6.}
$$
\begin{array}{|l|l|l|l|l|l|l|l|}
\hline
\alpha^3 &0	&1	&2	&3	&4	&5	&6 \\
\hline
0&	3&	0&	2&	5&	1&	4&	6  \\
1&	6&	3&	1&	2&	5&	0&	4  \\
2&	0&	5&	4&	6&	3&	1&	2  \\
3&	4&	1&	3&	0&	2&	6&	5  \\
4&	2&	4&	0&	1&	6&	5&	3  \\
5&	5&	2&	6&	4&	0&	3&	1  \\
6&	1&	6&	5&	3&	4&	2&	0  \\
\hline
\end{array}
$$

 \centerline{Table 7.}
$$
\begin{array}{|l|l|l|l|l|l|l|l|}
\hline
\beta^{\,6} &0	&1	&2	&3	&4	&5	&6 \\
\hline
0&	2&	5&	3&	0&	1&	4&	6 \\
1&	1&	2&	6&	3&	5&	0&	4 \\
2&	4&	6&	0&	5&	3&	1&	2 \\
3&	3&	0&	4&	1&	2&	6&	5 \\
4&	0&	1&	2&	4&	6&	5&	3 \\
5&	6&	4&	5&	2&	0&	3&	1 \\
6&	5&	3&	1&	6&	4&	2&	0 \\
\hline
\end{array}
$$

   \centerline{Table 8.}
$$
\begin{array}{|l|l|l|l|l|l|l|l|}
\hline
(\gamma^{5})^{(-1)} &0	&1	&2	&3	&4	&5	&6 \\
\hline
0&	6&	1&	0&	4&	3&	2&	5  \\
1&	3&	6&	5&	0&	1&	4&	2  \\
2&	2&	5&	4&	1&	0&	3&	6  \\
3&	0&	4&	2&	3&	6&	5&	1  \\
4&	4&	3&	6&	2&	5&	1&	0  \\
5&	5&	2&	1&	6&	4&	0&	3  \\
6&	1&	0&	3&	5&	2&	6&	4  \\
\hline
\end{array}
$$

Encryption. To send a message $b=630512403$, Bob computes from the known  $T=(\alpha,\beta,\gamma)$:
$\alpha =(2 3 4)  (0 5 1 6)$;  $\beta=(0 3 2 1)  (5 6)$; $\gamma=(1 2 3 6 0 5 4)$,
  calculates isotopy $T^{(r,s,t) }$  for random numbers $r=5, s=3, t=6$, i.e. $T^{(5,3,6)}$:

$
\alpha^5 =  \begin{pmatrix}
 0 &1& 2& 3& 4& 5& 6 \\
 5 &6 &4 &2 &3 &1 &0 \\
 \end{pmatrix}
$

$
\beta^3 =  \begin{pmatrix}
 0 &1& 2& 3& 4& 5& 6 \\
 1 &2 &3 &0 &4 &6 &5 \\
 \end{pmatrix}
$

$
\gamma^6 =  \begin{pmatrix}
 0 &1& 2& 3& 4& 5& 6 \\
 6 &4 &1 &2 &5 &0 &3 \\
 \end{pmatrix}
$

In our example  $T^{(5,3,6)}$ we get:
$\alpha^5=(0 5 1 6)  (2 4 3) $;   $\beta^3=(0 1 2 3)  (5 6)$;   $\gamma^6=(0 6 3 2 1 4 5)$.

Then he calculates $T^{(mr,ns,kt)}$ using the public key:

$T^{(m,n,k)}= (\alpha^m,\beta^n,\gamma^k)= (\alpha^*,\beta^*,\gamma^*)$:

$\alpha^*=\begin{pmatrix}
 0 &1& 2& 3& 4& 5& 6 \\
6& 5& 2& 3& 4& 0& 1 \\
 \end{pmatrix}
$

$\beta^*=\begin{pmatrix}
 0 &1& 2& 3& 4& 5& 6 \\
2 &3 &0 &1 &4 &5 &6 \\
 \end{pmatrix}
$
$\gamma^*=\begin{pmatrix}
 0 &1& 2& 3& 4& 5& 6 \\
3 &5 &4 &1 &0 &6 &2\\
 \end{pmatrix}
$

Then he raises these permutations, respectively, to the power $r=5,  s=3, t=6$ and gets:

$\alpha^{*5}=\begin{pmatrix}
 0 &1& 2& 3& 4& 5& 6 \\
6& 5& 2& 3& 4& 0& 1 \\
 \end{pmatrix}
$
$\beta^{*3}=\begin{pmatrix}
 0 &1& 2& 3& 4& 5& 6 \\
2 &3 &0 &1 &4 &5 &6 \\
 \end{pmatrix}
$

$\gamma^{*6}=\begin{pmatrix}
 0 &1& 2& 3& 4& 5& 6 \\
4 &3 &6 &0 &2 &1 &5\\
 \end{pmatrix}
$

$\alpha^{5m}=\alpha^5=(0 6 1 5)$;   $\beta^{3n}=\beta^3=(0 2)  (1 3)$; $  \gamma^{6k}=(0 4 2 6 5 1 3)$,
$(\gamma^{6k})^{(-1)}=(0 3 1 5 6 2 4)$.

As a result of the application of the new isotopy $T^{(5m,3n,6k)}$ to the quasigroup $(Q,f)$ we obtain:

   \centerline{Table 9.}
$$
\begin{array}{|l|l|l|l|l|l|l|l|}
\hline
\alpha^{5m} &0	&1	&2	&3	&4	&5	&6 \\
\hline
0&	3&	0&	2&	5&	1&	4&	6 \\
1&	6&	3&	1&	2&	5&	0&	4 \\
2&	0&	5&	4&	6&	3&	1&	2 \\
3&	4&	1&	3&	0&	2&	6&	5 \\
4&	2&	4&	0&	1&	6&	5&	3 \\
5&	5&	2&	6&	4&	0&	3&	1 \\
6&	1&	6&	5&	3&	4&	2&	0 \\
\hline
\end{array}
$$

  \centerline{Table 10.}
$$
\begin{array}{|l|l|l|l|l|l|l|l|}
\hline
\beta^{3n} &0	&1	&2	&3	&4	&5	&6 \\
\hline
0&	2&	5&	3&	0&	1&	4&	6 \\
1&	1&	2&	6&	3&	5&	0&	4 \\
2&	4&	6&	0&	5&	3&	1&	2 \\
3&	3&	0&	4&	1&	2&	6&	5 \\
4&	0&	1&	2&	4&	6&	5&	3 \\
5&	6&	4&	5&	2&	0&	3&	1 \\
6&	5&	3&	1&	6&	4&	2&	0 \\
\hline
\end{array}
$$

  \newpage

\centerline{Table 11.}
$$
\begin{array}{|l|l|l|l|l|l|l|l|}
\hline
(\gamma^{6k})^{(-1)} &0	&1	&2	&3	&4	&5	&6 \\
\hline
0&	4&	6&	1&	3&	5&	0&	2 \\
1&	5&	4&	2&	1&	6&	3&	0 \\
2&	0&	2&	3&	6&	1&	5&	4 \\
3&	1&	3&	0&	5&	4&	2&	6 \\
4&	3&	5&	4&	0&	2&	6&	1 \\
5&	2&	0&	6&	4&	3&	1&	5 \\
6&	6&	1&	5&	2&	0&	4&	3 \\
\hline
\end{array}
$$

To obtain $(T^{(mr,ns,kt)} (Q,f)b$, Bob uses the Markovski algorithm known to Alice, with the known leader value $l=3$, then the ciphertext for $b = 6\,3\,0\,5\,1\,2\,4\,0\,3$ will look like:
$v_1=3\cdot 6=6,
v_2= 6\cdot 3=2,
v_3= 2\cdot 0=0,
v_4= 0\cdot 5=0,
v_5= 0\cdot 1=6,
v_6= 6\cdot 2=5,
v_7=5 \cdot 4=3,
v_8= 3\cdot 0=1,
v_9= 1\cdot 3=1,  b^{\prime}=620065311$.

                                           \medskip

Decryption. Alice knows $m = 3,n = 6,k = 5$, so if she gets an isotopy  $T^{(r,s,t)}$  and ciphertext $(T^{(mr,ns,kt)} (Q,f))b)=620065311$,
she will calculate the isotopy first $T^{(mr,ns,kt)}$  using $T^{(r,s,t)} =T^{(**,**,**)}$ :

$\alpha^{\ast\ast} =\begin{pmatrix}
 0 &1& 2& 3& 4& 5& 6 \\
5 &6 &4& 2& 3& 1& 0 \\
 \end{pmatrix}
$

$\beta^{\ast\ast} =\begin{pmatrix}
 0 &1& 2& 3& 4& 5& 6 \\
1 &2 &3 &0 &4 &6 &5 \\
 \end{pmatrix}
$

$\gamma^{\ast\ast} =\begin{pmatrix}
 0 &1& 2& 3& 4& 5& 6 \\
6& 4& 1& 2& 5& 0& 3 \\
 \end{pmatrix}
$

She calculates $T^{(mr,ns,kt)}$:

$\alpha^{**3}=\begin{pmatrix}
 0 &1& 2& 3& 4& 5& 6 \\
6& 5& 2& 3& 4& 0& 1 \\
 \end{pmatrix}
$

$\beta^{**6}=\begin{pmatrix}
 0 &1& 2& 3& 4& 5& 6 \\
2 &3 &0 &1 &4 &5 &6 \\
 \end{pmatrix}
$

$\gamma^{**5}=\begin{pmatrix}
 0 &1& 2& 3& 4& 5& 6 \\
4 &3 &6 &0 &2 &1 &5 \\
 \end{pmatrix}
$

As a result, she receives the same table as Bob received in the encryption process.
For Table $(\gamma^{6k})^{-1} $ Alice builds a parastrophe $(23)$ used in the Markovski algorithm for decryption:
\centerline{Table 12.}
$$
\begin{array}{|l|l|l|l|l|l|l|l|}
\hline
\backslash &0	&1	&2	&3	&4	&5	&6 \\
\hline
0&	5&	2&	6&	3&	0&	4&	1 \\
1&	6&	3&	2&	5&	1&	0&	4 \\
2&	0&	4&	1&	2&	6&	5&	3 \\
3&	2&	0&	5&	1&	4&	3&	6 \\
4&	3&	6&	4&	0&	2&	1&	5 \\
5&	1&	5&	0&	4&	3&	6&	2 \\
6&	4&	1&	3&	6&	5&	2&	0 \\
\hline
\end{array}
$$
 and finally, using this table, she calculates b:
$u_1=3\backslash  6=6,
u_2= 6\backslash 2=3,
u_3= 2\backslash 0=0,
u_4= 0\backslash 0=5,
u_5= 0\backslash 6=1,
u_6= 6\backslash 5=2,
u_7=5 \backslash 3=4,
u_8= 3\backslash 1=0,
u_9= 1\backslash 1=3.$

Therefore, $b =630512403$.

\end{example}

In this algorithm, isostrophy [6] can also be used instead of isotopy, the modified algorithm instead of the Markovski algorithm and n-ary ($n> 2$) quasigroups [7; 8] instead of binary quasigroups.

A generalization of the Diffie-Hellman scheme of the open key distribution is given in [9].

The generalization is based on the concepts of the left and right powers of the elements of some non-associative groupoids.

For medial quasigroups, this approach is implemented in [10]. The protocol of the elaboration of a common secret key based on Moufang loops is given in [10].

This protocol is a generalization of the results from [11]. Generalizations of the ElGamal scheme based on Moufang loops are given in [10].

In [12], the discrete logarithmic problem with Moufang loops is reduced to the same problem over finite simple fields. Another generalization of the ElGamal scheme based on quasi-automorphisms of quasigroups is presented in [10].

\section{Conclusion}

Today, different points of view on the same mathematical idea lead to different generalizations. We considered in our work an analogue of the ElGamal encryption system based on the Markovski algorithm. This algorithm is under improvement and its other modifications are planned.

\noindent Nadeghda Malyutina$^{1}$

 $^{1}$  Ph.D. Student/Moldova State University

\noindent E--mail: \texttt{231003.bab.nadezhda@mail.ru}

\vspace{3mm}

\noindent Victor  Shcherbacov$^{2}$

$^{2}$ Researcher/Vladimir Andrunachievici Institute of Mathematics and Computer Science

\noindent
E--mail: \texttt{victor.scerbacov@math.md}

\end{document}